\documentclass[12pt,reqno]{amsart}
\usepackage{amsthm,amsmath,amssymb,amsfonts,amscd,verbatim,latexsym}

\newtheorem{Theorem}{Theorem}[section]
\newtheorem{Lemma}[Theorem]{Lemma}

\renewcommand{\P}{\mathbb P}
\newcommand{\E}{\mathcal E}
\newcommand{\I}{\mathbb I}
\newcommand{\J}{\mathbb J}
\newcommand{\tJ}{\widetilde{\mathbb J}}
\newcommand{\tQ}{\widetilde{\mathbb Q}}
\newcommand{\B}{\mathbb B}
\newcommand{\Q}{\mathbb Q}
\newcommand{\A}{\mathcal A}
\newcommand{\T}{\mathbb T}
\newcommand{\G}{{\mathbf G}^\circ}
\newcommand{\complex}{\mathbf C}
\newcommand{\F}{\mathbb F}
\newcommand{\ra}{\rightarrow}
\newcommand{\la}{\leftarrow}
\newcommand{\lra}{\longrightarrow}
\newcommand{\Sym}{\text{Sym}}
\newcommand{\Proof}{\noindent {\sc Proof.} \;} 
\newcommand{\key}{($\sharp$)}
\begin{document} 
\title[Ideals of binary orbits]{On the ideals of general binary orbits: \\ The low order cases} 
\author[Jaydeep Chipalkatti]{Jaydeep Chipalkatti${}^\star$}
\setcounter{footnote}{-1}
\footnote{\hspace{-4mm} ${}^\star$Department of Mathematics, 
University of Manitoba, Winnipeg, MB R3T 2N2, Canada. 
e-mail: chipalka@cc.umanitoba.ca}

\maketitle 

\parbox{11.8cm}{\small 
{\sc Abstract.} 
Let $E$ denote a general complex binary form of order $d$ (seen as a point in $\P^d$ ), 
and let $\Omega_E \subseteq \P^d$ denote the closure of its $SL_2$-orbit. In this note, we calculate 
the equivariant minimal generators of its defining ideal 
$I_E \subseteq \complex[a_0,\dots,a_d]$ for $4 \leqslant d \leqslant 10$. 
In order to effect the calculation, 
we introduce a notion called the `graded threshold character' of $d$. One unexpected feature of the 
problem is the (rare) occurrence of the so-called `invisible' generators in the ideal, and the 
resulting dichotomy on the set of integers $d \geqslant 4$.}

\thispagestyle{empty}

\bigskip \medskip 

\parbox{12cm}{\small 
} 

\bigskip  

\setcounter{footnote}{1}
\section{Introduction}
\subsection{} 
Let 
\[ E = \sum\limits_{i=0}^d \, \binom{d}{i} \, \alpha_i \, x_1^{d-i} \, x_2^i, \qquad (\alpha_i \in \complex) \] 
denote a nonzero form of order $d$ in the variables $\{x_1, x_2\}$. We will identify 
$E$ (distinguished up to a scalar) with the point 
$[\alpha_0, \dots, \alpha_d]$ in $\P^d$. 
Define the graded polynomial ring $R = \complex[a_0,\dots,a_d]$ over indeterminates 
$a_i$, so that $\P^d = \text{Proj} \, R$. 

The special linear group $SL_2  \complex$ acts on $R$ (and hence on $\P^d$) as follows. 
Let 
\begin{equation} \F = \sum\limits_{i=0}^d \, \binom{d}{i} \, a_i \, x_1^{d-i} \, x_2^i \, , 
\label{generic.dform} \end{equation} 
denote the generic binary $d$-ic. Given 
$g = \left( \begin{array}{cc} p & q \\ r & s \end{array} \right) \in SL_2$, make substitutions 
\[ x_1 = p \, x_1' + q \, x_2', \quad x_2 = r \, x_1' + s \, x_2',  \] 
into the right-hand side of~(\ref{generic.dform}), and rearrange terms to write 
\begin{equation} 
\F = \sum\limits_{i=0}^d   \, \binom{d}{i} \, a_i' \, {x_1'}^{d-i} \, {x_2'}^i. 
\label{genform.F} \end{equation} 
Then the action of $g$ takes $a_i$ to $a_i'$. 

\subsection{} 
If $d \geqslant 4$, and $E$ is a general point in $\P^d$, then the closure of its $SL_2$-orbit 
(denoted $\Omega_E$) is an irreducible projective variety of dimension $3$. The degree of $\Omega_E$ is 
$6$ for $d=4$, and $d \, (d-1)(d-2)$ for $d >4$ (see~\cite[p.~206]{AF} or~\cite[\S 8]{MJ}). 
Its defining ideal $I_E$ is an $SL_2$-subrepresentation of $R$, and we should like to find 
the equivariant minimal set of generators for $I_E$. 
This is similar (but not identical) to the `equivalence problem for binary forms.' (For 
discussions of the latter, see~\cite[\S 92]{Clebsch} or~\cite[Chapter 8]{Olver}.) 

The object of this paper is a complete determination of such generators for orders 
$d \leqslant 10$. The results are phrased in the language of classical invariant theory, i.e., in terms of 
invariants and covariants of the generic form $\F$. 

The Betti numbers of $I_E$ can be calculated by straightforward elimination 
(implemented here in Macaulay-2); it is rather the identification of the Betti modules 
{\sl qua} $SL_2$-representations which accounts for the bulk of 
the effort. In order to accomplish this, we introduce a notion 
called the {\sl graded threshold character} of $d$. Broadly speaking, it is designed to encode those 
subrepresentations of the ideal which can be detected by purely combinatorial considerations. 
This allows us to deduce an inequality involving the representation-theoretic character of a Betti module. 
It is a very surprising circumstance (to the author) that it turns out to be 
an {\sl equality} sufficiently often for the calculation to succeed. 

\section{Preliminaries} 
The ansatz used in this paper is similar to the one in~\cite[\S1]{JC1}, and the reader will find there 
detailed explanations of many of the notions used below. We refer to~\cite[Lecture 11]{FH} and 
~\cite[Chapter~4]{Sturmfels} for the basic representation theory of $SL_2$. Classical accounts of the 
invariant theory of binary forms may be found in~\cite{GY,Salmon1}, and more modern expositions 
in~\cite{Dolgachev,Olver,Springer}.  For the necessary facts from commutative 
algebra, reference~\cite{Eisenbud} is more than adequate. 

\subsection{} 
The base field will be $\complex$. Let $S_q$ denote the $(q+1)$-dimensional vector space of 
binary forms of order $q$ in $\{x_1,x_2\}$. Then $\{S_q: q \geqslant 0\}$ is the totality of all finite dimensional 
irreducible $SL_2$-representations. Since $SL_2$ is a linearly reductive group, each finite dimensional 
representation decomposes as a direct sum of irreducibles. We will need two specific 
decomposition formulae: the Clebsch-Gordan formula 
\begin{equation} 
S_p \otimes S_q \simeq \bigoplus\limits_{r=0}^{\min(p,q)} \, S_{p+q-2r}, 
\label{CG.formula} \end{equation}
and the Cayley-Sylvester formula 
\begin{equation} 
\Sym^p(S_q) \simeq \, \bigoplus\limits_{r=0}^{[\frac{p q}{2}]} \; 
(S_{pq-2r})^{\pi(r,p,q)-\pi(r-1,p,q)} \, . 
\label{CS.formula} \end{equation}
Here $\pi(a,b,c)$ denotes the number of partitions of $a$ into $b$ parts such that 
no part exceeds $c$. 

\subsection{} 
Given forms $A \in S_p$ and $B \in S_q$, the image of $A \otimes B$ via the projection map 
$S_p \otimes S_q \lra S_{p+q-2r}$ is called their $r$-th transvectant, denoted by $(A,B)_r$. 
It is given by the formula 
\[ (A,B)_r = \frac{(p-r)! \, (q-r)!}{p! \, q!} \, \sum\limits_{i=0}^r \, (-1)^i \, \binom{r}{i} \, 
\frac{\partial^r A}{\partial x_1^{r-i} \, \partial x_2^i} \, 
\frac{\partial^r B}{\partial x_1^i \, \partial x_2^{r-i}}. \] 
Two forms $A,B \in S_p$ are said to be {\sl apolar} if $(A,B)_p=0$. The pairing 
$S_p \otimes S_p \lra S_0$ is nondegenerate, hence there is a $p$-dimensional space of 
forms apolar to any specific nonzero form $A \in S_p$ (see~\cite[Chapter XI]{GY}). 
\subsection{} 
Let $\Gamma$ denote the representation ring of $SL_2$, i.e., it is a free abelian group on generators 
$s_0,s_1,s_2,\dots$ etc.,~with multiplication corresponding to the tensor product of representations.  
Given a finite-dimensional $SL_2$-representation $U$, let 
$[U] \in \Gamma$ denote its character. For instance, 
\[ [S_5 \otimes S_3] = s_5 \cdot s_3 = s_8 + s_6 + s_4 + s_2, \] 
by formula~(\ref{CG.formula}). We will write $s_p \circ s_q$ for $[\Sym^p(S_q)]$, e.g., 
\begin{equation} \begin{aligned} 
s_4 \circ s_7 = \, & s_{28} + s_{24} + s_{22} + 2 \, s_{20} + s_{18} + 3 \, s_{16} + 
 2 \, s_{14} + \\ & 3 \, s_{12} + 2 \, s_{10} + 3 \, s_8 + s_6 + 3 \, s_4 + s_0, 
\end{aligned} \label{s4os7} \end{equation} 
by formula~(\ref{CS.formula}). 

Given elements 
$a = \sum\limits_i \alpha_i \, s_i$ and $b = \sum\limits_i  \beta_i \, s_i$ in $\Gamma$, 
write $a \geqslant b$ if $\alpha_i \geqslant \beta_i$ for all $i$. Define 
\[ \sup(a,b) = \sum\limits_i \, \max(\alpha_i,\beta_i) \, s_i. \] 

\subsection{} 
There is an isomorphism of $SL_2$-representations 
\[ \jmath: R_1 \stackrel{\sim}{\ra} S_d, \quad a_i \ra (-1)^i \, x_1^i \, x_2^{d-i};  \] 
which allows us to make the identification 
\[ R = \bigoplus\limits_{m \geqslant 0} \, R_m = \bigoplus\limits_{m \geqslant 0} \, \Sym^m(S_d). \]  
Let $W_{m,q} \subseteq R_m$ denote the span of the images of all $SL_2$-equivariant maps 
$S_q \lra R_m$. Then there is a decomposition of representations 
\[ R_m = \bigoplus\limits_q \;  W_{m,q}. \] 
Similarly $(I_E)_m = \bigoplus\limits_q \; (I_E)_{m,q}$, where 
$(I_E)_{m,q}$ is an $SL_2$-invariant subspace of $W_{m,q}$. 

\subsection{} 
Let $\A_{m,q}$ denote the space of covariants of $\F$ in degree $m$ and 
order $q$. Each element $\Phi \in \A_{m,q}$ may be written as  
\[ \sum\limits_{i=0}^q \, \varphi_i(a_0,\dots,a_d) \, x_1^{q-i} \, x_2^i, \] 
where $\varphi_i$ are homogeneous forms of degree $m$ in $a_0,\dots,a_d$. 
Now $\Phi$ defines an equivariant morphism 
\[ S_q \lra R_m, \quad A(x_1,x_2)  \lra (\Phi,A)_q, \] 
whose image is $\text{Span} \, \{\varphi_0,\dots,\varphi_q\} \subseteq W_{m,q}$. Every such 
morphism comes from a covariant, i.e., we have an isomorphism 
\[ \A_{m,q} \simeq \text{Hom}_{SL_2}(S_q, W_{m,q}).  \] 
This induces a bijection between subspaces of $\A_{m,q}$ and $SL_2$-invariant subspaces of $W_{m,q}$. 
It associates to a subspace $U \subseteq \A_{m,q}$, the span of all the coefficients of all the 
elements in $U$ (to be denoted by $U^\circ$). 

It is a standard fact (see~\cite[\S 86]{GY}) that $\A_{m,q}$ admits a basis each of whose 
elements is a compound transvectant in $\F$. E.g., for $d=7$, the space $\A_{3,9}$ is $2$-dimensional 
with a basis $\{\F \, (\F,\F)_6, (\F,(\F,\F)_2)_4\}$. By formula~(\ref{CS.formula}), 
\[ \zeta_{m,q} = \dim \A_{m,q} = \pi(\frac{m \, d-q}{2}, m, d) - \pi(\frac{m \, d-q-2}{2}, m, d). \] 

\subsection{} Given a specific form $E \in S_d$, there is an evaluation map 
\[ \theta_E: \A_{m,q} \lra S_q \] 
which substitutes the coefficients of $E$ for the indeterminates $a_i$. 
Write $(K_E)_{m,q} = \ker \theta_E$. Henceforth we may omit $E$ from the notation if no confusion is likely; 
it is understood that $K,I,J$ etc depend upon the choice of $E$. 

\begin{Lemma} \sl 
We have an equality $K_{m,q}^\circ = I_{m,q}$. 
\end{Lemma} 
\Proof Indeed, each element in $K_{m,q}^\circ$ vanishes on $E$, and hence 
by equivariance on $\Omega_E$. Alternately, let $e \in I_{m,q}$ denote a nonzero element. 
Then $e$ belongs to a unique smallest $SL_2$-invariant subspace $V \subseteq I_{m,q}$. 
Let $\Psi \in \A_{m,q}$ denote the covariant (unique up to a constant) whose coefficients give a 
basis of $V$. It immediately follows that $\Psi \in K_{m,q}$, hence $e \in K_{m,q}^\circ$. \qed 

\medskip 

Since $\dim K_{m,q}$ is no smaller than $\max (0,\zeta_{m,q}-q-1)$, 
we will define the {\sl threshold character (of $d$) in degree} $m$ to be the element
\[ \T_m = \sum\limits_{q \geqslant 0} \; \max \, (0,\zeta_{m,q}-q-1) \, s_q \in \Gamma. \] 

For instance, let $d=5$. Then $\zeta_{14,10} = 17$, hence the coefficient of 
$s_{10}$ in $\T_{14}$ is $17-11=6$. In fact, the full expression is 
\begin{equation} 
 \T_{14} = s_{22} + 4 \, s_{18} + 2 \, s_{16} + 6 \, s_{14} + 3 \, s_{12} + 
6 \, s_{10} + 2 \, s_8 + 5 \, s_6 + s_4 + 3 \, s_2. \label{T5.14} \end{equation}

\subsection{} 
The minimal resolution of $I$ will be written as 
\[ 0 \la I \la \E_0 \la \E_1 \la \cdots \] 
where 
$\E_r = \bigoplus\limits_{j \geqslant 0} \, B(r,j) \otimes R(-j)$, 
are graded free $R$-modules of finite rank. Thus $B(0,-)$ are the minimal generators of $I_E$, and 
$B(r,-)$ are the $r$-th syzygy modules. Each Betti module $B(r,j)$ is 
an $SL_2$-representation. For each $m$, let $J_m \subseteq I_m$ denote the subspace generated by 
$\bigoplus\limits_{r < m} I_r$ (the ideal elements in earlier degrees). Write 
\[ \I_m = [I_m], \quad \J_m = [J_m], \quad \B_m = [B(0,m)] \] 
for the corresponding elements in $\Gamma$. By construction, $\I_m \geqslant \T_m$ (which justifies the 
term `threshold'), and hence 
\[ \B_m = \I_m - \J_m \geqslant \underbrace{\sup(\J_m, \T_m) - \J_m}_{\Q_m}. \] 
Henceforth $E$ will be assumed to be sufficiently general, which ensures 
that $\I,\J,\B$ etc are independent of $E$. 
\subsection{} The Betti numbers (in the free resolution of) $I$ can be calculated as follows. 
For illustration, let $d=6$. Choose a `general' form $E(x_1,x_2)$ in $S_6$, and make 
simultaneous substitutions 
\[ x_1 \ra p \, x_1 + q \, x_2, \quad x_2 \ra r \, x_1 + s \, x_2, \] 
into $E$ to construct a new form 
\[ \sum\limits_{i=0}^6 \, \binom{6}{i} \, \psi_i(p,q,r,s) \, x_1^{6-i} \, x_2^i. \] 
This defines a ring morphism 
\[ \Psi_E: \complex[a_0,\dots,a_6] \lra \complex[p,q,r,s], \quad a_i \lra \psi_i(p,q,r,s). \] 
Then $I= \ker \Psi_E$. The actual calculation shows that the Betti numbers of $I$ are 
as in the following table: 
\[ \begin{array}{r|cccccc}
{} & 0 & 1 & 2 & 3 & 4 & 5 \\ \hline 
4   & 1 \\ 
6   & 1 \\ 
9  &      & 1 \\ 
10 & 1 \\
12 & 97 & 222 & 114 &  7  & 1  \\ 
13 & 27 & 235 & 609 & 587 & 233 & 30 \\ 
15 & & & & & & 1 
\end{array} \] 
The entry in the row labelled $i$ and column labelled $j$ gives the dimension of 
$B(j,i+j)$, e.g., $\dim B(1,14) = 235$. In practice, for each $d$, I have repeated the 
calculation for several random choices of $E$ to eliminate any likelihood of error. 
Our task is to identify the $B(0,m)$ {\sl qua} $SL_2$-representations, and 
secondly to identify the corresponding ideal generators. 

\subsection{} It is a paradoxical feature of the subsequent calculations that the higher 
syzygies do not enter into them.\footnote{With one small exception, noted later.} Define 
\[ \tJ_m = [(\E_0)_m] - \B_m = \sum\limits_{j < m}  \;  [ B(0,j) \otimes R_{m-j}] = 
\sum\limits_{j < m}  \; \B_j \cdot (s_{m-j} \circ s_d).  \] 
This should be thought of as an approximation to $\J_m$, but with all higher syzygies ignored. 
Clearly $\tJ_m \geqslant \J_m$. Define $\tQ_m = \sup(\tJ_m,\T_m) - \tJ_m$. 
\begin{Lemma} \sl We have an inequality $\Q_m \geqslant \tQ_m$. 
\end{Lemma} 
\Proof Fix an integer $q$, and let $a,b,c$ denote the coefficients of $s_q$ in 
$\tJ_m, \J_m$ and $\T_m$ respectively. Since $a \geqslant b$, the result follows from 
the obvious inequality $\max(b,c)-b \geqslant \max(a,c)-a$. \qed 

\medskip 

As a consequence, we have the crucial inequality 
\[ \B_m \geqslant \tQ_m   \qquad (\sharp) \] 
which will serve as our workhorse throughout the next section. 
\section{Computations} 
In this section we will describe the solution for each $d$. The calculations for order $d$ are 
to be found in \S 3.d. Of course, the results are valid only for $E$ belonging to a dense open subset of $\P^d$. E.g., 
if $E = x_1^d$, then $\Omega_E$ is the rational normal curve whose ideal is generated by quadrics. 

Henceforth we will write $\beta_m$ for $\dim B(0,m)$, to be called the generator dimensions of $I$. 
As mentioned earlier, they were all calculated using {\sc Macaulay-2}. 
Formulae~(\ref{CG.formula}), (\ref{CS.formula}) as well as the rest of the 
calculations in the representation ring $\Gamma$ were programmed in {\sc Maple} by the author. 

We will determine the $\B_m$ successively for increasing $m$. If the characters 
$\B_r$ for $r<m$ are known, then the calculation of $\tQ_m$ is a purely mechanical task. 
Now our governing principle is simple: if the dimensions of $\B_m$ and 
$\tQ_m$ coincide\footnote{The dimension of an element in $\Gamma$ is understood in the obvious sense.}, 
then we must have equality in~\key. At first blush, this seems optimistic beyond reason. 
However, it is an intriguing but pleasing circumstance that~\key~is an equality in all the cases 
below, with only two exceptions. Moreover, each of the exceptions is `thematic' in a sense which will 
be readily understood once it is encountered. 

We will say that {\sl all the ideal generators in degree $m$ are visible} if~\key~is an 
equality; if not, the ideal $I$ is said to have {\sl invisible generators in degree $m$}. These phrases are 
to be understood atomically; it is meaningless to speak of any specific element in the ideal 
as being visible or otherwise. 

\setcounter{subsection}{3}

\subsection{Quartics} \label{section.quartic}
The variety $\Omega_E$ is a hypersurface of degree $6$. Since $\zeta_{6,0} = 2$, 
the space $K_{6,0}$ is one-dimensional, and its generator gives the defining equation for $\Omega_E$. 
Said differently, define invariants 
\[ g_2 = (\F,\F)_4, \quad g_3 = (\F,(\F,\F)_2)_4, \] 
in degrees $2$ and $3$ respectively. Then $\{g_2^3,g_3^2\}$ is a basis of $\A_{6,0}$, and hence 
\[ \theta_E(g_3)^2 \, g_2^3 - \theta_E(g_2)^3 \, g_3^2 \] is the generator of $I$. 
\subsection{Quintics} \label{section.quintic}
The generator dimensions are 
\[ \beta_8 =1, \quad \beta_{12} = 1, \quad \beta_{14} = 60. \] 
Evidently, $\B_8 = \B_{12} = s_0$. Define covariants (cf.~\cite[p.~131]{GY}) 
\[ \begin{array}{lll} 
H = (\F,\F)_2, &  \imath = (\F,\F)_4, \\ 
A = (\imath,\imath)_2, & B = (\imath^3, H)_6, & C = (\imath^5, \F^2)_{10}. 
\end{array} \] 
Now $\zeta_{8,0}=2, \, \zeta_{12,0}=3$, and 
\[ \{A^2, B \}, \quad \{A^3, A \, B, C\} \] 
are respectively bases of $\A_{8,0}$ and $\A_{12,0}$. As in the previous section, the degree $8$ 
generator of $I$ can be taken to be $Z_8 = \theta_E(B) \, A^2 - \theta_E(A)^2 \, B$. Then the 
new generator in degree $12$ can be chosen to be any element in $K_{12,0}$ which is not a 
multiple of $A \, Z_8$, e.g., $Z'_{12} = \theta_E(A \, B) \, C - \theta_E(C) \, A \, B$.  

Now, $\tJ_{14} = s_6 \circ s_5 + s_2 \circ s_5$, which evaluates to 
\[ \begin{aligned} 
s_{30} & + s_{26} + s_{24} + 2 \, s_{22} + 2 \, s_{20} + 3 \, s_{18} + 2 \, s_{16} + 
4 \, s_{14} + 3 \, s_{12} \, + \\ 5 \, s_{10} & + 2 \, s_8 + 5 \, s_6 +  s_4 + 3 \, s_2, 
\end{aligned} \] 
by formula~(\ref{CS.formula}). Using the expression for $\T_{14}$ from~(\ref{T5.14}), one arrives at 
\[ \tQ_{14} = s_{18} + 2 \, s_{14} + s_{10}, \] 
which has dimension $19 + 2 \cdot 15 + 11 = 60$. Hence $\B_{14} 
= \tQ_{14}$ in \key. 

We introduce some notation in order to describe the generators succintly. 
There is a $2$-dimensional subspace $V \subseteq K_{14,14}$, such that $V^\circ$ accounts for the 
new generators in degree $14$. Now $V$ is not uniquely determined, but all choices satisfying the 
condition $V^\circ \cap J_{14} = (0)$ are valid. Henceforth we will write $\G(2,K_{14,14})$ 
for such a $V^\circ$. 
In general, $\G(r,K)$ will stand for the span of coefficients of 
an $r$-dimensional subspace $V \subseteq K$, where $V$ is chosen to lie outside 
a (tacitly specified) proper subvariety in the Grassmannian of $r$-subspaces of $K$. 

We have arrived at the following result: 

\begin{Theorem} \sl 
For a general quintic $E$, the ideal $I$ is minimally generated by the following subspaces: 
\[ K_{8,0}, \quad \G(1,K_{12,0}), \quad \G(1,K_{14,18}), \quad 
\G(2,K_{14,14}), \quad \G(1,K_{14,10}). \] 
\end{Theorem} 

\subsection{Sextics} \label{section.sextic}
The generator dimensions are 
\[ \begin{array}{lllll} 
\beta_4 = 1, & \beta_6 = 1, & \beta_{10} = 1, & 
\beta_{12} = 97, & \beta_{13} = 27. \end{array} \] 
Hence $\B_4 = \B_6 = \B_{10} =s_0$. A preliminary man{\oe}uvre is necessary before 
proceeding to degree $12$. Notice that the generators in degrees $4,6$ must give rise to a first syzygy 
in degree $10$. Its contribution to $\I_m$ can be cancelled against that of the degree $10$ generator. 
Thus, for the purposes of calculating $\tJ_m$ and $\tQ_m$, we will henceforth ignore $\B_{10}$. 
Then one gets
\[ \tQ_{12} = s_{24} + 2 \, s_{20} + s_{16} + s_{12}, \] 
which is $97$-dimensional; hence $\B_{12} = \tQ_{12}$. Similarly, 
\[ \B_{13} = \tQ_{13} = s_{26},  \] 
which completes the calculation. We have proved the following result. 
\begin{Theorem} \sl 
For a general sextic $E$, the ideal $I$ is minimally generated by the following subspaces: 
\[ \begin{array}{llll} 
K_{4,0},  & \G(1,K_{6,0}), & \G(1,K_{10,0}), & \G(1,K_{12,24}), \\ 
\G(2,K_{12,20}), & \G(1,K_{12,16}), & \G(1,K_{12,12}), & \G(1,K_{13,26}). 
\end{array} \] 
\end{Theorem}                   
We will no longer state such theorems explicitly, since they can be written down ritually once 
the $\B_m$ are known. 

\subsection{Septimics} \label{section.septimic}
The generator dimensions are 
\[ \beta_6 = 10, \quad \beta_8 = 40, \quad \beta_9 = 106, \quad \beta_{10} = 89. \] 
 
A calculation shows that $\T_6$ (and hence $\tQ_6$) equals $s_2$. It follows 
that~\key \, must be a strict inequality, i.e., there are invisible generators in degree $6$. 
The explanation lies behind the following algebraic peculiarity 
of the ring of covariants for binary septimics. 

The spaces $\A_{4,6}$ and $\A_{6,6}$ are respectively of dimensions $1$ and $7$. Let 
$\Delta$ denote a generator of the former.\footnote{One can choose 
$\Delta = ((\F,\F)_4, (\F,\F)_6)_1$, but the precise expression is not relevant to the argument.}
Septimics have no invariant in degree $10$, i.e., $\A_{10,0}=0$. 
It follows that for any $\Phi \in \A_{6,6}$, we must have $(\Phi,\Delta)_6=0$. But then 
$(\theta_E(\Phi),\theta_E(\Delta))_6 = 0$, i.e., the image of the evaluation map 
\[ \theta_E: \A_{6,6} \lra S_6 \] 
is contained in the $6$-dimensional subspace of sextics which are apolar to 
$\theta_E(\Delta)$. Hence $K_{6,6} \neq 0$. It follows that $s_6$ must be a summand 
in $\B_6$, and hence on dimensional grounds $\B_6 = s_6 + s_2$. 

The rest of the generators are all visible, hence the calculation is straightforward. The Betti 
modules are 
\[ \begin{aligned} 
{} & \B_8=s_{16} + s_{12} + s_8 + s_0, \\ 
& \B_9=s_{23} + 2 \, s_{21} + s_{19} + s_{17}, \\ 
& \B_{10} = 2 \, s_{30} + s_{26}. 
\end{aligned} \] 

\subsection{Octavics}  \label{section.octavic}
The generator dimensions are 
\[ \begin{array}{llll} 
\beta_4 = 1, & \beta_5 = 1, & \beta_6 = 7, & \beta_7 = 106, \\ 
\beta_8 = 264, & \beta_9 = 97,  & \beta_{10} = 82. \end{array} \] 
All the generators are visible, and the Betti modules are 
\[ \begin{aligned} 
{} & \B_4= s_0, \\ 
& \B_5 = s_0, \\ 
& \B_6=s_4+2 \, s_0, \\ 
& \B_7=2 \, s_{16} + s_{14} + 2 \, s_{12} + s_{10} + s_8 + 2 \, s_4 + s_0, \\
& \B_8=4 \, s_{24} + 2 \, s_{22} + 4 \, s_{20} + 2 \, s_{16}, \\
& \B_9=2 \, s_{32} + s_{30}, \\ 
& \B_{10} = 2 \, s_{40}. \end{aligned} \] 
\subsection{Nonics} \label{section.nonic}
The generator dimensions are 
\[ \begin{array}{llll} 
\beta_4 = 1, & \beta_6 = 71, & \beta_7 = 508, & \beta_8 = 324, \\ 
\beta_9 = 86, & \beta_{10} = 51. \end{array} \] 
Once again, all the generators are visible. The Betti modules are 
\[ \begin{aligned} 
{} & \B_4=s_0, \\ 
& \B_6=s_{14}+2 \, s_{10} + 4 \, s_6 + 2 \, s_2, \\ 
& \B_7=3 \, s_{23} + 5 \, s_{21} + 5 \, s_{19} + 6 \, s_{17} + 4 \, s_{15} 
+ 3 \, s_{13} + s_{11}. \\ 
& \B_8=s_{34} + 6 \, s_{32} + 2 \, s_{30} + s_{28}, \\ 
& \B_9=s_{43} + s_{41}, \\ 
& \B_{10} =s_{50}. \end{aligned} \]

\subsection{Decimics}  \label{section.decimic}
The generator dimensions are 
\[ \begin{array}{llll}
\beta_4 = 1, & \beta_5 = 3, & \beta_6 = 367, & \beta_7 = 679,  \\ 
\beta_8 = 324, & \beta_9 = 151, & \beta_{10} = 61. 
\end{array} \] 
The generators in degrees $4,5$ are visible, which gives 
$\B_4 = s_0$ and $\B_5 = s_2$. In degree $6$, one gets 
\[ \tQ_6 = 
2 \, s_{20} + 5 \, s_{16} + 2 \, s_{14} + 7 \, s_{12} + s_{10} + 6 \, s_8 + 2 \, s_6 + 
5 \, s_4 + 4 \, s_0,  \] 
which is $356$-dimensional, hence there exist invisible generators. The explanation is similar to the case of 
septimics. 

The space $\A_{6,10}$ is $13$-dimensional, whereas $\A_{7,0}=0$. 
Thus every element in $\A_{6,10}$ is apolar to $\F$. It follows that the map 
$\A_{6,10} \stackrel{\theta_E}{\lra} S_{10}$ is not surjective, and hence its kernel 
is at least $3$-dimensional. The coefficient of 
$s_{10}$ in $\tJ_6$ is $1$, hence $\B_6$ must contain at least two copies of $s_{10}$. This forces 
$\B_6 = \tQ_6 + s_{10}$, since the additional term precisely compensates for the missing 
dimensions ($367=356+11$). 

From degree $7$ onwards, all the generators are visible and the modules are 
\[ \begin{aligned} 
{} & \B_7=4 \, s_{30} + 6 \, s_{28} + 7 \, s_{26} + 4 \, s_{24} + 4 \, s_{22}, \\ 
& \B_8=6 \, s_{40} + 2 \, s_{38}, \\ 
& \B_9=2 \, s_{50} + s_{48}, \\ 
& \B_{10}=s_{60}. \end{aligned} \] 

\medskip 

I know of no general method for identifying the characters corresponding to 
invisible generators. In either of the cases above, 
it is only by educated guesswork that we have succeeded in doing so. 
 
\section{Miscellaneous Remarks} 

\subsection{} 
Let us say (for the present purposes) that an integer $d \geqslant 4$ is `prosaic' if \key~is an 
equality for all $m$, and `erratic' otherwise. Our calculations show that $d=4,5,8,9$ are prosaic, 
whereas $d=7,10$ are erratic.

We have treated the case $d=6,m=10$ as anomalous. Following the 
definition literally, one gets $\tQ_{10}=0$, i.e., we have strict inequality in \key. Nevertheless, 
(as we have seen) it is easy to restore equality by cancelling $\B_{10}$ against a first syzygy. 
This suggests that our definitions of `prosaic' and `erratic' are not in their final shape, and 
a more refined understanding of the problem will modify them. However, 
even in their present formulation they do seem to capture a valuable distinction. 

It would be an interesting (but immensely ambitious) undertaking to \nobreak{arrive} at 
such a classification for all $d$. The problem implicitly involves the structure of the 
ring of covariants $\bigoplus\limits_{m,q} \A_{m,q}$. Such rings are in general very 
complicated, and it is not obvious how to proceed in the general case. 

\subsection{} The process we have used to calculate the ideal generators is analogous to 
the minimal resolution conjecture ({\sc mrc}) for general points in $\P^n$ (see~\cite{Lorenzini}). 
To see the parallel, consider the following example: let $X$ denote a set of $8$ general points in 
$\P^2$, and we are to find the generator degrees of its defining 
ideal $I_X \subseteq R = \complex[z_0,z_1,z_2]$. The heuristic reasoning goes as follows. 
Since $\dim R_3 = 10$, the evaluation map 
$e_X: R_3 \lra \complex^8$ has kernel dimension $\geqslant 2$. Since the points are general, 
we may assume equality, i.e., $\dim \, (I_X)_3 = 2$. 
By the same reasoning, $\dim \, (I_X)_4 = 15-8 = 7$. Now one assumes that 
the rank of the map $(I_X)_3 \otimes R_1 \lra (I_X)_4$ is the maximum possible, which is $2 \times 3 = 6$. 
Hence there should be one new generator in degree $4$. The process detects no further 
generators in degree $5$, hence we have an expected presentation 
\[ 0 \la R/I_X  \la R \la R(-3)^2 \oplus R(-4) \la \ldots \] 

The argument can be continued to obtain the module of first syzgygies of $I_X$ (which would be 
$R(-5)^2$ in this case), but I have not succeeded in the analogous calculation for $I_E$. Although 
{\sc mrc} is false in general (see~\cite{EP}), it is known to be true in many cases 
(in particular for $\P^2$). Thus, broadly speaking, the dichotomy between prosaic and erratic integers 
corresponds to the one between true and false instances of {\sc mrc}. 

\subsection{} 
There is an evidently analogous problem of calculating $I_E$ 
for the action of $SL_n$ on the space of $n$-ary $d$-ics. To the best 
of my knowledge, the answer is known only in the case $d=n=3$. Ternary cubics have two 
invariants $G_4, G_6$  in degrees $4,6$ respectively (cf.~\cite[\S 198]{Salmon2}, where 
they are labelled $S$ and $T$). 
For a general cubic curve $E$, the hypersurface 
$\Omega_E \subseteq \P^9$ is of degree $12$, with defining equation 
\[ \theta_E(G_6)^2 \, G_4^3 - \theta_E(G_4)^3 \, G_6^2 = 0. \] 
Much to my chagrin, I have found that at present even the case of ternary quartics 
seems too large for computational experimentation. 

\subsection{} 
If one considers the same problem (in the binary case) over a field of 
characteristic $p >0$, then preliminary calculations show that the generator dimensions 
of $I_E$ depend on $p$. Here are some data for $d=5$. 
\[ \begin{array}{c|l} 
\text{characteristic}  & \\ \hline 
2 & \beta_8=\beta_{12}=1, \, \beta_{13}=\beta_{14}=12, \, \beta_{16}=18 \\ 
3 & \beta_8=\beta_{12}=1, \, \beta_{13}=6, \, \beta_{14}=32, \, \beta_{15}=6 \\ 
5 & \text{same as characteristic zero} \\ 
7 & \beta_8=\beta_{12}=1, \, \beta_{13}=2, \, \beta_{14}=48 \\ 
11 & \text{same as characteristic zero} \\
13 & \text{same as characteristic zero} 
\end{array} \] 

Since $SL_2$ is no longer linearly reductive, many of the techniques used here are no longer applicable.

\bigskip 

{\small \noindent {\sc acknowledgements:} The author is thankful to Daniel Grayson and Michael 
Stillman (the authors of Macaulay-2). This work was financially supported by a discovery grant from 
NSERC, Canada.}

\medskip \centerline{---} 
\end{document}